\documentclass[reqno,a4paper]{amsart}%
\usepackage{amssymb,url}%,lineno
\usepackage[colorlinks=true,hypertex]{hyperref}%,pagebackref
\date{Completed on 1 October 2008 at Student Village of Victoria University in Australia}
\date{}
\allowdisplaybreaks[4]
\theoremstyle{plain}
\newtheorem{thm}{Theorem}
\theoremstyle{definition}
\newtheorem{dfn}{Definition}
\theoremstyle{remark}
\newtheorem{rem}{Remark}
\DeclareMathOperator{\td}{d\mspace{-1mu}}
\begin{document}

\title[The function $(b^x-a^x)/x$: Logarithmic convexity and applications]
{The function $\boldsymbol{(b^x-a^x)/x}$: Logarithmic convexity and applications to extended mean values}

\author[F. Qi]{Feng Qi}
\address[F. Qi]{Department of Mathematics, College of Science, Tianjin Polytechnic University, Tianjin City, 300160, China}
\email{\href{mailto: F. Qi <qifeng618@gmail.com>}{qifeng618@gmail.com}, \href{mailto: F. Qi <qifeng618@hotmail.com>}{qifeng618@hotmail.com}, \href{mailto: F. Qi <qifeng618@qq.com>}{qifeng618@qq.com}}
\urladdr{\url{http://qifeng618.spaces.live.com}}

\author[B.-N. Guo]{Bai-Ni Guo}
\address[B.-N. Guo]{School of Mathematics and Informatics, Henan Polytechnic University, Jiaozuo City, Henan Province, 454010, China}
\email{\href{mailto: B.-N. Guo <bai.ni.guo@gmail.com>}{bai.ni.guo@gmail.com}, \href{mailto: B.-N. Guo <bai.ni.guo@hotmail.com>}{bai.ni.guo@hotmail.com}}
\urladdr{\url{http://guobaini.spaces.live.com}}

\begin{abstract}
In the present paper, we first prove the logarithmic convexity of the elementary function $\frac{b^x-a^x}x$, where $x\ne0$ and $b>a>0$. Basing on this, we then provide a simple proof for Schur-convex properties of the extended mean values, and, finally, discover some convexity related to the extended mean values.
\end{abstract}



\thanks{The first author was supported partially by the China Scholarship Council}

\thanks{This paper was typeset using \AmS-\LaTeX}

\maketitle

\section{Introduction}

For given numbers $b>a>0$, let
\begin{equation}
g_{a,b}(t)=
\begin{cases}
\dfrac{b^t-a^t}t,&t\not=0;\\
b-a,&t=0.
\end{cases}
\end{equation}
This elementary and special function was first dedicated to be investigated in~\cite{jmaa-ii-97, (b^x-a^x)/x}. Subsequently, it was utilized to construct Steffensen pairs in~\cite{Gauchman-Steffensen-pairs, steffensen-qi-cheng-rgmia, steffensen-pair-Anal, qcw, onsp} and its reciprocal was also used to generalize Bernoulli numbers and polynomials in \cite{ijmest-bernoulli, bernoulli-luo-guo-qi-rgmia, bernoulli-luo-guo-qi-debnath-IJMMS, euler-bernoulli-luo-qi-adv, bernoulli-qi-guo-rgmia}. It has something to do with the classical Euler gamma function $\Gamma$ and the remainder of Binet's first formula for the logarithm of $\Gamma$ (see, for example, \cite{best-constant-one-simple.tex, remiander-Sen-Lin-Guo.tex, best-constant-one.tex, Extended-Binet-remiander-comp.tex, Extended-Binet-remiander-comp.tex-Slovaca, notes-best.tex, psi-reminders.tex-ijpams, psi-reminders.tex-rgmia, best-constant-one-simple-real.tex} and closely-related references therein). More importantly, it was employed not only to provide alternative proofs for the monotonicity of the extended mean values $E(r,s;x,y)$ in \cite{Qi-Luo-Simple-98, NEW-L.TEX} but also to create the logarithmic convexity and Schur-convex properties of $E(r,s;x,y)$ in \cite{Cheung-Qi-mean-rgmia, Cheung-Qi-mean, emv-rs-schur-rocky, pams-62-rgmia, pams-62, emv-rs-schur-rgmia, Cheung-Qi-Rev.tex}, where the extended mean values $E(r,s;x,y)$ were defined in \cite{ls1, stolarsky} for $x,y>0$ and $r,s\in\mathbb{R}$ as
\begin{align*}
&E(r,s;x,y)=\biggl(\frac{r}{s}\cdot\frac{y^s-x^s} { y^r-x^r}\biggr)^{{1/(s-r)}}, & & rs(r-s)(x-y)\ne 0; \\
&E(r,0;x,y)=\biggl(\frac{1}{r}\cdot\frac{y^r-x^r}{ \ln y-\ln x}\biggr)^{{1/r}}, & & r(x-y)\ne 0; \\
&E(r,r;x,y)=\frac1{e^{1/r}}\biggl(\frac{x^{x^r}}{y^{y^r}}\biggr)^{ {1/(x^r-y^r)}},& & r(x-y)\ne 0; \\
&E(0,0;x,y)=\sqrt{xy}, & & x\ne y; \\
&E(r,s;x,x)=x, & & x=y.
\end{align*}
There has been a lot of literature on the extended mean values $E(r,s;x,y)$. For more information, please refer to \cite{bullen-handbook, cubo} and related references therein.
\par
In this paper, we first present the logarithmic convexity of the function $g_{a,b}(t)$. Basing on this, we then provide a concise proof for Schur-convex properties of the extended men values $E(r,s;x,y)$, and, finally, discover some monotonicity and logarithmic convexity of certain functions related to the extended mean values $E(r,s;x,y)$.

\section{Logarithmic convexity of $g_{a,b}(t)$}

For the sake of proceeding smoothly, we need the following definition which can be found in~\cite{note-on-neuman.tex, note-on-neuman-ITSF-simplified.tex} and related references therein.

\begin{dfn}
A $k$-times differentiable function $f(t)>0$ is said to be $k$\nobreakdash-log\nobreakdash-convex on an interval $I$ if
\begin{equation}\label{k-times-ineq}
[\ln f(t)]^{(k)}\ge0,\quad k\in\mathbb{N}
\end{equation}
on $I$; If the inequality~\eqref{k-times-ineq} reverses then $f$ is said to be $k$\nobreakdash-log\nobreakdash-concave on $I$.
\end{dfn}

Now we are in a position to state and prove the logarithmic convexity of the function $g_{a,b}(t)$ on $(-\infty,\infty)$.

\begin{thm}\label{exp-f-f-T2}
Let $b>a>0$. Then the function $g_{a,b}(t)$ is logarithmic convex on $(-\infty,\infty)$, $3$-log-convex on $(-\infty,0)$, and $3$-log-concave on $(0,\infty)$. Consequently, the function
\begin{equation}
h_{a,b}(t)=
\begin{cases}
\dfrac{b^t\ln b-a^t\ln a}{b^t-a^t}-\dfrac1t,&t\ne0\\
\ln\sqrt{ab}\,,&t=0
\end{cases}
\end{equation}
is increasing on $(-\infty,\infty)$ and satisfies
\begin{equation}
\lim_{t\to-\infty}h_{a,b}(t)=\ln a\quad \text{and} \quad \lim_{t\to\infty}h_{a,b}(t)=\ln b.
\end{equation}
\end{thm}

\begin{proof}
For $t\ne0$, taking the logarithm of $g_{a,b}(t)$ and differentiating yields
\begin{align*}
\ln g_{a,b}(t)&=\ln\bigl\vert b^t-a^t\bigr\vert-\ln\vert{t}\vert,\\
[\ln g_{a,b}(t)]'&=\frac{b^t\ln b-a^t\ln a}{b^t-a^t}-\frac1t,\\
[\ln g_{a,b}(t)]''&=\frac1{t^2}-\frac{a^tb^t(\ln a-\ln b)^2}{(a^t-b^t)^2}
\end{align*}
and
\begin{multline*}
[\ln g_{a,b}(t)]'''=\frac{a^tb^t(a^t+b^t)(\ln a-\ln b)^3} {(a^t-b^t)^3}-\frac2{t^3}\\
\begin{aligned}
&=\frac{2(ab)^{3t/2}}{t^3}\biggl(\frac{t\ln a-t\ln b}{a^t-b^t}\biggr)^3 \biggl\{\frac{(a/b)^{t/2}+(b/a)^{t/2}}2 -\biggl[\frac{(a/b)^{t/2}-(b/a)^{t/2}}{(\ln a-\ln b)t}\biggr]^3\biggr\}\\
&\triangleq \frac{2(ab)^{3t/2}}{t^3}\biggl(\frac{t\ln a-t\ln b}{a^t-b^t}\biggr)^3 Q_{a,b}(t),
\end{aligned}
\end{multline*}
where, by using Lazarevi\'c's inequality in \cite[p.~131]{bullen-dic} and \cite[p.~300]{3rded},
\begin{equation*}
Q_{a,b}\biggr(\frac{2t}{\ln a-\ln b}\biggr) =\frac{e^{-t}+e^t}{2}-\biggl(\frac{e^t-e^{-t}}{2t}\biggr)^3 =\cosh t-\biggl(\frac{\sinh t}{t}\biggr)^3<0.
\end{equation*}
Consequently,
\begin{equation*}
[\ln g_{a,b}(t)]'''=\biggl[\frac{g_{a,b}'(t)}{g_{a,b}(t)}\biggr]''
\begin{cases}
>0, & t\in(-\infty,0)\\
<0, & t\in(0,\infty)
\end{cases}
\end{equation*}
which implies that the function $[\ln g_{a,b}(t)]''$ is increasing on $(-\infty,0)$ and decreasing on $(0,\infty)$. Since
\begin{equation*}
\lim_{t\to-\infty}\frac{a^tb^t}{(a^t-b^t)^2} =\lim_{t\to-\infty}\frac{(a/b)^t}{[(a/b)^t-1]^2} =\lim_{t\to-\infty}\frac{(b/a)^t}{[(b/a)^t-1]^2}=0
\end{equation*}
and the function $[\ln g_{a,b}(t)]''$ is even on $\mathbb{R}$, then $[\ln g_{a,b}(t)]''>0$, and so the function $[\ln g_{a,b}(t)]'=h_{a,b}(t)$ is increasing on $\mathbb{R}$. Since
\begin{equation*}
\frac{b^t\ln b-a^t\ln a}{b^t-a^t} =\frac{(b/a)^t\ln b-\ln a}{(b/a)^t-1} =\frac{\ln b-(a/b)^t\ln a}{1-(a/b)^t},
\end{equation*}
then it follows easily that
\begin{equation*}
\begin{aligned}
\lim_{t\to-\infty}[\ln g_{a,b}(t)]'&=\ln a & &\text{and} &
\lim_{t\to\infty}[\ln g_{a,b}(t)]'&=\ln b.
\end{aligned}
\end{equation*}
The L'H\^ospital's rule reveals that
\begin{equation*}
\begin{split}
\lim_{t\to0}\{[\ln g_{a,b}(t)]'\}
&=\lim_{t\to0}\frac{t(b^t\ln b-a^t\ln a)-(b^t-a^t)}{t(b^t-a^t)}\\
&=\lim_{t\to0}\frac{y^t(\ln b)^2-x^t(\ln a)^2}{(b^t-a^t)/t +(b^t\ln b-a^t\ln a)}\\
&=\frac{\ln b+\ln a}2.
\end{split}
\end{equation*}
The proof of Theorem~\ref{exp-f-f-T2} is thus completed.
\end{proof}

\begin{rem}
In the preprint \cite{exp-funct-further.tex}, Theorem~\ref{exp-f-f-T2} was also verified by using the celebrated Hermite-Hadamard's integral inequality \cite{difference-hermite-hadamard.tex, correction-to-sandor.tex-ijams, correction-to-sandor.tex-octogon, correction-to-sandor.tex-rgmia} instead of Lazarevi\'c's inequality.
\end{rem}

\begin{rem}
Theorem~\ref{exp-f-f-T2} provides important supplements to the work in~\cite{jmaa-ii-97, (b^x-a^x)/x}.
\end{rem}

\section{A concise proof of Schur-convexity of $E(r,s;x,y)$}

Let us recall \cite[pp.~75\nobreakdash--76]{cposa} the definition of Schur-convex functions.

\begin{dfn}\label{schur-dfn}
A function $f$ with $n$ arguments defined on $I^n$ is called Schur\nobreakdash-convex if $f(x)\le f(y)$ holds for each two $n$-tuples $x=(x_1,\dotsc,x_n)$ and $y=(y_1,\dotsc,y_n)$ on $I^n$ such that $x \prec y$, where $I$ is an interval with nonempty interior and the relationship of majorization $x\prec y$ means that
\begin{equation}
\begin{aligned}
\sum_{i=1}^kx_{[i]}&\le\sum_{i=1}^ky_{[i]}& &\text{and}&
\sum_{i=1}^nx_{[i]}&=\sum_{i=1}^ny_{[i]}
\end{aligned}
\end{equation}
for $1\le k\le n-1$, where $x_{[i]}$ denotes the $i$-th largest component in $x$.
\par
A function $f$ is Schur\nobreakdash-concave if and only if $-f$ is Schur\nobreakdash-convex.
\end{dfn}

Based on intricate conclusions in \cite{pams-62, pams-62-rgmia} and basic properties of $E(r,s;x,y)$, the following Schur\nobreakdash-convex properties of the extended mean values $E(r,s;x,y)$ with respect to $(r,s)$ was first obtained in \cite{emv-rs-schur-rocky, emv-rs-schur-rgmia}.

\begin{thm}\label{logschur1}
With respect to the $2$-tuple $(r,s)$, the extended mean values $E(r,s;x,y)$ are Schur\nobreakdash-concave on $[0,\infty)\times[0,\infty)$ and Schur\nobreakdash-convex on $(-\infty,0]\times(-\infty,0]$.
\end{thm}

The aim of this section is to demonstrate a concise proof of Theorem~\ref{logschur1} with the help of Theorem~\ref{exp-f-f-T2}.

\begin{proof}
When $y>x>0$, the extended mean values $E(r,s;x,y)$ may be represented in terms of $g_{x,y}(t)$ as
\begin{equation*}
E(r,s;x,y)=
\begin{cases}
\label{case1}\biggl[\dfrac{g_{x,y}(s)}{g_{x,y}(r)}\biggr]^{1/(s-r)}, & (r-s)(x-y)\neq 0; \\[1em]
\exp \biggl[\dfrac{g_{x,y}'(r)}{g_{x,y}(r)}\biggr], & r=s,\; x-y\neq 0
\end{cases}
\end{equation*}
and
\begin{equation}\label{1.9}
\ln E(r,s;x,y)=
\begin{cases}
\dfrac{1}{s-r}\displaystyle\int_{r}^{s}\dfrac{g_{x,y}'(u)}{g_{x,y}(u)}\td u, & (r-s)(x-y)\neq 0; \\[1em]
\dfrac{g_{x,y}'(r)}{g_{x,y}(r)}, & r=s,\; x-y\neq 0.
\end{cases}
\end{equation}
In virtue of Theorem~\ref{exp-f-f-T2}, it follows that
\begin{equation}\label{ln-g-3-der-sign}
[\ln g_{x,y}(t)]^{(3)}=\biggl[\frac{g_{x,y}'(t)}{g_{x,y}(t)}\biggr]''
\begin{cases}
<0, &t\in(0,\infty),\\
>0, &t\in(-\infty,0).
\end{cases}
\end{equation}
\par
In \cite{Elezovic-Pacaric-rocky-2000}, it was obtained that the integral arithmetic mean
\begin{equation}  \label{lem11}
\phi(r,s)= \begin{cases}
\dfrac1{s-r}\displaystyle\int_r^sf(t)\td t, & r\ne s \\
f(r), & r=s
\end{cases}
\end{equation}
of a continuous function $f$ on $I$ is Schur\nobreakdash-convex (or Schur\nobreakdash-concave, respectively) on $I^2$ if and only if $f$ is convex (or concave, respectively) on $I$. Consequently, by virtue of the formula~\eqref{1.9} and Definition~\ref{schur-dfn}, it is not difficult to see that, in order that the extended mean values $E(r,s;x,y)$ are Schur-convex (or Schur-concave, respectively) with respect to $(r,s)$, it is sufficient to show the validity of \eqref{ln-g-3-der-sign}, which may be deduced from Theorem~\ref{exp-f-f-T2} straightforwardly. Theorem~\ref{logschur1} is thus proved.
\end{proof}

\begin{rem}
In \cite{sandor-banach-07}, an alternative proof of Theorem~\ref{logschur1} was given, among other things.
\end{rem}

\section{Some logarithmic convexity related to $E(r,s;x,y)$}

In Remark~6 of~\cite{new-upper-kershaw.tex, new-upper-kershaw-JCAM.tex}, it was pointed out that the reciprocal of the exponential mean
\begin{equation}
I_{s,t}(x)=\frac1e\biggl[\frac{(x+s)^{x+s}}{(x+t)^{x+t}}\biggr]^{ {1/(s-t)}}
\end{equation}
for $s\ne t$ is logarithmically completely monotonic on $(-\min\{s,t\},\infty)$ and that the exponential mean $I_{s,t}(x)$ for $s\ne t$ is also a completely monotonic function of first order on $(-\min\{s,t\},\infty)$.
\par
In \cite{gamma-psi-batir.tex}, it was remarked that the logarithmic mean
\begin{equation}
L_{s,t}(x)=L(x+s,x+t)
\end{equation}
is increasing and concave on $(-\min\{s,t\},\infty)$ for $s\ne t$. In \cite{log-mean-comp-mon.tex, log-mean-comp-mon.tex-mia}, the logarithmic mean $L_{s,t}(x)$ for $s\ne t$ is further proved to be a completely monotonic function of first order on $(-\min\{s,t\},\infty)$.
\par
For $x,y>0$ and $r,s\in\mathbb{R}$, let
\begin{align}
F_{r,s;x,y}(w)&=E(r+w,s+w;x,y),\quad w\in\mathbb{R},\\
G_{r,s;x,y}(w)&=E(r,s;x+w,y+w),\quad w>-\min\{x,y\}
\end{align}
and
\begin{equation}
H_{r,s;x,y}(w)=E(r+w,s+w;x+w,y+w),\quad w>-\min\{x,y\}.
\end{equation}
By virtue of the monotonicity of the extended mean values $E(r,s;x,y)$, it is easy to see that the functions $F_{r,s;x,y}(w)$, $G_{r,s;x,y}(w)$ and $H_{r,s;x,y}(w)$ are increasing with respect to $w$. Furthermore, since
$$
I_{s,t}(x)=E(1,1;x+s,y+t)=G_{1,1;x,y}(w)
$$
and
$$
L_{s,t}(x)=E(0,1;x+s,y+t)=G_{0,1;x,y}(w),
$$
the following problem was posed in \cite{log-mean-comp-mon.tex, log-mean-comp-mon.tex-mia}: What about the logarithmic convexity of the functions $F_{r,s;x,y}(w)$, $G_{r,s;x,y}(w)$ and $H_{r,s;x,y}(w)$ with respect to $w$?
\par
The aim of this section is to supply a solution to the above problem about the function $F_{r,s;x,y}(w)$. Our main results are the following theorems.

\begin{thm}\label{log-comp-mon-log-thm}
The function $F_{r,s;x,y}(w)$ is logarithmically convex on $\bigl(-\infty,-\frac{s+r}2\bigr)$ and logarithmically concave on $\bigl(-\frac{s+r}2,\infty\bigr)$.
\end{thm}

\begin{thm}\label{log-comp-mon-log-thm-2}
The product $\mathcal{F}_{r,s;x,y}(w)=F_{r,s;x,y}(w)F_{r,s;x,y}(-w)$ is increasing on $(-\infty,0)$ and decreasing on $(0,\infty)$.
\end{thm}

\begin{thm}\label{Cheung-Qi-Rev-conv}
If $s+r>0$, the function $w\ln F_{r,s;x,y}(w)$ is convex on $\bigl(-\frac{s+r}2,0\bigr)$; If $s+r<0$, it is also convex on $\bigl(0,-\frac{s+r}2\bigr)$.
\end{thm}

\begin{proof}[Proof of Theorem~\ref{log-comp-mon-log-thm}]
In the first place, we claim that if $f(t)$ is even on $(-\infty,\infty)$ and increasing on $(-\infty,0)$, then the function
\begin{equation}\label{f-f=F}
p(t)=f(t+\alpha)-f(t),\quad \alpha>0
\end{equation}
is positive on $\bigl(-\infty,-\frac\alpha2\bigr)$ and negative on $\bigl(-\frac\alpha2,\infty\bigr)$. This can be verified as follows:
\begin{enumerate}
\item
If $t+\alpha>t>0$, since $f(t)$ is decreasing on $(0,\infty)$, then $F(t)<0$;
\item
If $t<t+\alpha<0$, since $f(t)$ is increasing on $(-\infty,0)$, then $F(t)>0$;
\item
If $t+\alpha>0>t$,
\begin{enumerate}
\item
when $t+\alpha>-t>0$, i.e., $t>-\frac\alpha2$, using the even and monotonic properties of $f(t)$ shows that $F(t)=f(t+\alpha)-f(-t)$ and it is negative;
\item
similarly, when $-t>t+\alpha>0$, i.e., $t<-\frac\alpha2$, the function $F(t)$ is positive.
\end{enumerate}
\end{enumerate}
The claim is thus proved.
\par
From~\eqref{1.9}, it follows that if $y>x>0$ then
\begin{equation}\label{emv-w}
\frac{\td^2\ln F_{r,s;x,y}(w)}{\td w^2}=
\begin{cases}\displaystyle
\dfrac{1}{s-r}\int_{r}^{s}\frac{\td^2}{\td w^2} \biggl[\frac{g_{x,y}'(w+t)}{g_{x,y}(w+t)}\biggr]\td t, &(r-s)(x-y)\ne0;\\[1em]
\dfrac{\td^2}{\td w^2} \biggl[\dfrac{g_{x,y}'(w+r)}{g_{x,y}(w+r)}\biggr], & r=s,\,\, x-y\ne 0.
\end{cases}
\end{equation}
As shown in the proof of Theorem~\ref{exp-f-f-T2}, the function $[\ln g_{x,y}(t)]''$ for $y>x>0$ is even on $\mathbb{R}$ and increasing on $(-\infty,0)$. Substituting $f(t)$ and $\alpha$ by $[\ln g_{x,y}(t)]''$ and $s-r>0$ in \eqref{f-f=F} respectively and utilizing~\eqref{emv-w} demonstrates that
\begin{equation*}
\frac{[\ln g_{x,y}(t+s-r)]''_t-[\ln g_{x,y}(t)]''}{s-r}
=\frac{\td^2\ln F_{r,s;x,y}(t-r)}{\td t^2}>0
\end{equation*}
for $t<-\frac{s-r}2$ and that $\frac{\td^2\ln F_{r,s;x,y}(t-r)}{\td t^2}<0$ for $t>-\frac{s-r}2$. As a result,
\begin{equation}\label{emv-w-positive}
\frac{\td^2\ln F_{r,s;x,y}(w)}{\td w^2}
=\frac{[\ln g_{x,y}(w+s)]''_w-[\ln g_{x,y}(w+r)]''_w}{s-r}
\begin{cases}
>0,&w<-\dfrac{s+r}2,\\[0.5em]
<0,&w>-\dfrac{s+r}2.
\end{cases}
\end{equation}
Because $F_{r,s;x,y}(w)=F_{r,s;y,x}(w)=F_{s,r;x,y}(w)$, the equation \eqref{emv-w-positive} holds for all $r,s\in\mathbb{N}$ and $x,y>0$ with $x\ne y$. Theorem~\ref{log-comp-mon-log-thm} is proved.
\end{proof}

\begin{proof}[Proof of Theorem~\ref{log-comp-mon-log-thm-2}]
It is easy to see that
\begin{equation*}
\bigl[\ln\mathcal{F}_{r,s;x,y}(w)\bigr]' =\frac{F_{r,s;x,y}'(w)}{F_{r,s;x,y}(w)} -\frac{F_{r,s;x,y}'(-w)}{F_{r,s;x,y}(-w)}.
\end{equation*}
Careful computation reveals that
\begin{align*}
\frac{F_{r,s;x,y}'(w)}{F_{r,s;x,y}(w)} =\frac{F_{r,s;x,y}'(-w-(s+r))}{F_{r,s;x,y}(-w-(s+r))}
\end{align*}
for $w\in(-\infty,\infty)$. Theorem~\ref{log-comp-mon-log-thm} implies that the function
\begin{equation*}
q(w)=\frac{F_{r,s;x,y}'(w-(s+r)/2)}{F_{r,s;x,y}(w-(s+r)/2)}
\end{equation*}
is increasing on $(-\infty,0)$ and decreasing on $(0,\infty)$. It is also apparent that the function $q(w)$ is even, that is, $q(w)=q(-w)$ for $w\in(-\infty,\infty)$. By virtue of the claim verified in the proof of Theorem~\ref{log-comp-mon-log-thm}, it is easy to see that the difference $q(w+(s+r))-q(w)$ is positive on $\bigl(-\infty,-\frac{s+r}2\bigr)$ and negative on $\bigl(-\frac{s+r}2,\infty\bigr)$, equivalently, the function
\begin{equation}\label{q-q-w}
q\biggl(w+\frac{s+r}2\biggr)-q\biggl(w-\frac{s+r}2\biggr) =\frac{F_{r,s;x,y}'(w)}{F_{r,s;x,y}(w)} -\frac{F_{r,s;x,y}'(w-(s+r))}{F_{r,s;x,y}(w-(s+r))}
\end{equation}
is positive on $(-\infty,0)$ and negative on $(0,\infty)$. On the other hand, since
\begin{equation}
\mathcal{F}_{r,s;x,y}(w)=\frac{xyF_{r,s;x,y}(w)}{F_{r,s;x,y}(w-(s+r))},
\end{equation}
then the function~\eqref{q-q-w} equals $\bigl[\ln\mathcal{F}_{r,s;x,y}(w)\bigr]'$. Thus, Theorem~\ref{log-comp-mon-log-thm-2} is proved.
\end{proof}

\begin{proof}[Proof of Theorem~\ref{Cheung-Qi-Rev-conv}]
Direct calculation yields
\begin{equation}\label{klm}
\bigl[w\ln F_{r,s;x,y}(w)\bigr]''=2\bigl[\ln F_{r,s;x,y}(w)\bigr]' +w\bigl[\ln F_{r,s;x,y}(w)\bigr]''.
\end{equation}
By Theorem~\ref{log-comp-mon-log-thm}, it follows that $\bigl[\ln F_{r,s;x,y}(w)\bigr]'>0$ on $(-\infty,\infty)$, $\bigl[\ln F_{r,s;x,y}(w)\bigr]''>0$ on $\bigl(-\infty,-\frac{s+r}2\bigr)$, and $\bigl[\ln F_{r,s;x,y}(w)\bigr]''<0$ on $\bigl(-\frac{s+r}2,\infty\bigr)$. Therefore,
\begin{enumerate}
\item
if $s+r<0$, then $\bigl[w\ln F_{r,s;x,y}(w)\bigr]''>0$, and so $w\ln F_{r,s;x,y}(w)$ is convex on $\bigl(0,-\frac{s+r}2\bigr)$;
\item
if $s+r>0$, then $\bigl[w\ln F_{r,s;x,y}(w)\bigr]''>0$, and so $w\ln F_{r,s;x,y}(w)$ is convex on $\bigl(-\frac{s+r}2,0\bigr)$.
\end{enumerate}
The proof of Theorem~\ref{Cheung-Qi-Rev-conv} is complete.
\end{proof}

\begin{rem}
Theorem~\ref{log-comp-mon-log-thm} generalizes \cite[Theorem~1 and Theorem~3]{Cheung-Qi-mean} and \cite[Theorem~1]{Cheung-Qi-Rev.tex}. Theorem~\ref{log-comp-mon-log-thm-2} generalizes \cite[Theorem~2 and Theorem~3]{Cheung-Qi-mean} and \cite[Theorem~2]{Cheung-Qi-Rev.tex}. Theorem~\ref{Cheung-Qi-Rev-conv} generalizes \cite[Theorem~5]{Cheung-Qi-mean}. Notice that the paper~\cite{Cheung-Qi-mean-rgmia} is a preprint and a complete version of~\cite{Cheung-Qi-mean}.
\end{rem}

\begin{rem}
By the same method as in \cite[Theorem~4]{Cheung-Qi-mean}, the function
\begin{equation}
(w+s-r)\bigl[F_{r,s;x,y}(w)\bigr]^{s-r},\quad s>r
\end{equation}
can be proved to be increasingly convex on $(-\infty,\infty)$ and logarithmically concave on $\bigl(-\frac{s-r}2,\infty\bigr)$.
\end{rem}

\end{document}